\documentclass[10pt,letterpaper]{article}
\usepackage[T1]{fontenc}
\usepackage{lmodern,microtype}
\usepackage[margin=0.6in,top=0.65in,bottom=0.65in]{geometry}
\usepackage{amsmath,amssymb,amsthm,graphicx,booktabs,multirow,array}
\usepackage[authoryear,round]{natbib}
\usepackage{xurl}
\usepackage[hidelinks,pdftitle={The Scaled Dirichlet Distribution: Some Properties},pdfauthor={Jared N. Lakhani}]{hyperref}
\allowdisplaybreaks[2]
\newtheorem{definition}{Definition}
\newtheorem{lemma}{Lemma}
\newtheorem{theorem}{Theorem}
\newtheorem{proposition}{Proposition}
\newtheorem{corollary}{Corollary}
\newcommand{\E}{\mathbb E}
\newcommand{\Cov}{\operatorname{Cov}}
\newcommand{\Var}{\operatorname{Var}}
\newcommand{\Cor}{\operatorname{Cor}}
\newcommand{\FD}{F_D}
\newcommand{\SD}{\operatorname{SDirichlet}}
\newcommand{\Dir}{\operatorname{Dirichlet}}
\newcommand{\Ga}{\operatorname{Gamma}}
\newcommand{\dd}{\,\mathrm d}
\newcommand{\ind}{\mathbf 1}
\newcommand{\bm}[1]{\boldsymbol{#1}}
\title{The Scaled Dirichlet Distribution: Lauricella Representations and Estimation}
\author{Jared N. Lakhani\\
\texttt{jaredlakhani@gmail.com}}
\date{}
\begin{document}
\maketitle
\begin{abstract}
We derive explicit expressions for the mixed moments of the scaled Dirichlet
distribution in terms of the Lauricella hypergeometric function of type~$D$.
The relationship with Carlson's multivariate $R$-function is stated, and
the convergence conditions required for the series derivation are made
explicit. The general expressions yield means, second moments and
covariances, with Gauss and Appell hypergeometric representations in the
bivariate and trivariate cases. Method-of-moments and maximum-likelihood formulations
are presented, together with simulation summaries and an illustrative
comparison of scaled Dirichlet and ordinary Dirichlet fits to compositional
futures-return data.
\end{abstract}

\section{The \texorpdfstring{$n$}{n}-Dimensional Scaled Dirichlet Distribution}
The purpose of this note is to collect explicit Lauricella expressions for
the scaled Dirichlet moments and their low-dimensional specialisations.
Their connection with established Dirichlet-average identities is recorded
in Corollary~\ref{cor:carlson} and Lemma~\ref{lem:carlson}. Related moment
representations are given by \citet{ng2011dirichlet} and \cite{Furman2019}.

\begin{definition}[Scaled Dirichlet Distribution; \cite{ng2011dirichlet}]
Let $n\geq2$ and
\[
 \bm x=(x_1,\ldots,x_n)^{\mathsf T}\sim\Dir_n(\bm a)
\]
on the simplex $T_n=\{\bm x:x_i>0,\ \sum_{i=1}^nx_i=1\}$. Define
\[
 w_i=\frac{b_ix_i}{\sum_{j=1}^n b_jx_j},\qquad i=1,\ldots,n,
\]
where $b_i>0$ for $i=1,\ldots,n$. Then the distribution of
$\bm w=(w_1,\ldots,w_n)^{\mathsf T}$ is called the \emph{scaled Dirichlet
distribution}, with positive shape parameter vector
$\bm a=(a_1,\ldots,a_n)^{\mathsf T}$ and scale parameter vector
$\bm b=(b_1,\ldots,b_n)^{\mathsf T}$. We write
$\bm w\sim\SD_n(\bm a,\bm b)$ on $T_n$.
\end{definition}
Since the components sum to one, only $n-1$ of them are free. Thus,
$T_n$ is an $(n-1)$-dimensional simplex. We write the density in terms
of $w_1,\ldots,w_{n-1}$, with $w_n=1-\sum_{i=1}^{n-1}w_i$:
\[
 f_{\bm W}(\bm w)=\frac{1}{B(\bm a)}
 \frac{\prod_{i=1}^n b_i^{-a_i}w_i^{a_i-1}}
 {\left(\sum_{j=1}^n w_j/b_j\right)^{a_+}},
\]
where
\[
 B(\bm a)=\frac{\prod_{i=1}^n\Gamma(a_i)}{\Gamma(a_+)},
 \qquad a_+=\sum_{i=1}^n a_i.
\]
The scale parameters must be positive, but they need not be less than
one. As shown in Section~2, multiplying all scales by the same positive
constant does not change the distribution.

\begin{lemma}[Denominator Representation]\label{lem:denominator}
Let
\[
 W_i=\frac{b_iX_i}{\sum_{j=1}^n b_jX_j},\qquad
 \bm X\sim\Dir_n(\bm a),
\]
and define $z_k=1-b_k/b_n$ for $k=1,\ldots,n-1$. Then
\begin{align*}
 \sum_{j=1}^n b_jX_j
 &=b_nX_n+\sum_{k=1}^{n-1}b_kX_k\\
 &=b_n\left(1-\sum_{k=1}^{n-1}X_k\right)
   +\sum_{k=1}^{n-1}b_kX_k\\
 &=b_n+\sum_{k=1}^{n-1}(b_k-b_n)X_k
 =b_n\left(1-\sum_{k=1}^{n-1}z_kX_k\right).
\end{align*}
Consequently,
\[
 W_i=\frac{b_i}{b_n}X_i\left(1-\sum_{k=1}^{n-1}z_kX_k\right)^{-1},
 \qquad i=1,\ldots,n.
\]
\end{lemma}

\begin{lemma}[Multivariate Binomial Series]\label{lem:binomial}
Let $p\in\mathbb N$ and let $u_1,\ldots,u_p$ satisfy
$\sum_{k=1}^p|u_k|<1$. For any $\lambda>0$,
\[
 \left(1-\sum_{k=1}^p u_k\right)^{-\lambda}
 =\sum_{m_1,\ldots,m_p\geq0}
 \frac{(\lambda)_{m_+}}{m_1!\cdots m_p!}\prod_{k=1}^p u_k^{m_k},
\]
where $m_+=m_1+\cdots+m_p$ and $(q)_m=q(q+1)\cdots(q+m-1)$ denotes
the Pochhammer symbol, with $(q)_0=1$.
\end{lemma}
\begin{proof}
Start from the ordinary binomial series
\[
 (1-u)^{-\lambda}=\sum_{m=0}^{\infty}\frac{(\lambda)_m}{m!}u^m,
 \qquad |u|<1.
\]
Set $u=\sum_{k=1}^p u_k$. Then
\[
 \left(1-\sum_{k=1}^p u_k\right)^{-\lambda}
 =\sum_{m=0}^{\infty}\frac{(\lambda)_m}{m!}
 \left(\sum_{k=1}^p u_k\right)^m.
\]
By the multinomial expansion in Lemma~\ref{lem:multinomial},
\[
 \left(\sum_{k=1}^p u_k\right)^m
 =\sum_{\substack{m_1,\ldots,m_p\geq0\\m_1+\cdots+m_p=m}}
 \frac{m!}{m_1!\cdots m_p!}\prod_{k=1}^p u_k^{m_k}.
\]
Substituting gives
\[
 \left(1-\sum_{k=1}^p u_k\right)^{-\lambda}
 =\sum_{m=0}^{\infty}\frac{(\lambda)_m}{m!}
 \sum_{\substack{m_1,\ldots,m_p\geq0\\m_1+\cdots+m_p=m}}
 \frac{m!}{m_1!\cdots m_p!}\prod_{k=1}^p u_k^{m_k}.
\]
Each term is determined by the tuple $(m_1,\ldots,m_p)$ and
$m=m_1+\cdots+m_p=m_+$. The assumption $\sum_k|u_k|<1$ gives absolute
convergence, so the sums may be combined into a single sum over all
nonnegative integer tuples. Cancelling $m!$ yields the result.
\end{proof}

\begin{theorem}[Mixed Moments of the Scaled Dirichlet Distribution]\label{thm:mixed}
Let $\bm W\sim\SD_n(\bm a,\bm b)$, with $a_i>0$, $b_i>0$, and let
$\bm r=(r_1,\ldots,r_n)$ with $r_i\in\mathbb N_0$. Define
$r_+=\sum_{i=1}^nr_i$, $a_+=\sum_{i=1}^na_i$ and
$z_k=1-b_k/b_n$ for $k=1,\ldots,n-1$. Then
\begin{equation}\tag{1}\label{eq:mixed_general}
 \boxed{\begin{aligned}
 \E\!\left[\prod_{i=1}^n W_i^{r_i}\right]
 &=\left[\prod_{i=1}^{n-1}\left(\frac{b_i}{b_n}\right)^{r_i}\right]
 \frac{\prod_{i=1}^n(a_i)_{r_i}}{(a_+)_{r_+}}\\[-2pt]
 &\quad\times\FD^{(n-1)}\!\left(r_+;
 a_1+r_1,\ldots,a_{n-1}+r_{n-1};a_++r_+;z_1,\ldots,z_{n-1}\right),
 \end{aligned}}
\end{equation}
where $\FD^{(m)}$ denotes the Lauricella hypergeometric function of type~$D$.
For $r_+>0$, its value is interpreted through the real Euler integral in
Definition~\ref{def:Lauricella} whenever the defining series is outside its
convergence domain. When $r_+=0$, both the moment and the Lauricella factor
equal one. Choosing the largest scale as the reference, after relabelling
if necessary, ensures $0\leq z_k<1$ and convergence of the defining series.
\end{theorem}
\begin{proof}
The case $r_+=0$ is immediate. First suppose that $b_n=\max_i b_i$.
By Lemma~\ref{lem:denominator},
\[
 \prod_{i=1}^nW_i^{r_i}
 =\left[\prod_{i=1}^{n-1}\left(\frac{b_i}{b_n}\right)^{r_i}\right]
 \left(\prod_{i=1}^nX_i^{r_i}\right)
 \left(1-\sum_{k=1}^{n-1}z_kX_k\right)^{-r_+}.
\]
Here $0\leq z_k<1$ and
$0\leq\sum_{k<n}z_kX_k\leq\max_{k<n}z_k<1$. Applying
Lemma~\ref{lem:binomial},
\[
 \left(1-\sum_{k=1}^{n-1}z_kX_k\right)^{-r_+}
 =\sum_{m_1,\ldots,m_{n-1}\geq0}
 \frac{(r_+)_{m_+}}{m_1!\cdots m_{n-1}!}
 \prod_{k=1}^{n-1}(z_kX_k)^{m_k},
\]
where $m_+=\sum_{k=1}^{n-1}m_k$. All terms are nonnegative, so Tonelli's
theorem justifies interchanging the expectation and the sum. Thus,
\begin{align*}
 \E\!\left[\prod_{i=1}^nW_i^{r_i}\right]
 &=\left[\prod_{i=1}^{n-1}\left(\frac{b_i}{b_n}\right)^{r_i}\right]
 \sum_{m_1,\ldots,m_{n-1}\geq0}
 \frac{(r_+)_{m_+}}{m_1!\cdots m_{n-1}!}
 \prod_{k=1}^{n-1}z_k^{m_k}
 \E\!\left[\prod_{k=1}^{n-1}X_k^{r_k+m_k}X_n^{r_n}\right].
\end{align*}
By Lemma~\ref{lem:moment_application},
\[
 \E\!\left[\prod_{k=1}^{n-1}X_k^{r_k+m_k}X_n^{r_n}\right]
 =\frac{\prod_{k=1}^{n-1}(a_k)_{r_k+m_k}(a_n)_{r_n}}
 {(a_+)_{r_++m_+}}.
\]
Using $(a_k)_{r_k+m_k}=(a_k)_{r_k}(a_k+r_k)_{m_k}$ and
$(a_+)_{r_++m_+}=(a_+)_{r_+}(a_++r_+)_{m_+}$ gives
\begin{align*}
 \E\!\left[\prod_{i=1}^nW_i^{r_i}\right]
 &=\left[\prod_{i=1}^{n-1}\left(\frac{b_i}{b_n}\right)^{r_i}\right]
 \frac{\prod_{i=1}^n(a_i)_{r_i}}{(a_+)_{r_+}}\\
 &\quad\times\sum_{m_1,\ldots,m_{n-1}\geq0}
 \frac{(r_+)_{m_+}\prod_{k=1}^{n-1}(a_k+r_k)_{m_k}}
 {(a_++r_+)_{m_+}}\prod_{k=1}^{n-1}\frac{z_k^{m_k}}{m_k!}.
\end{align*}
The remaining series is the Lauricella function in
Definition~\ref{def:Lauricella}, with first parameter $r_+$, component
parameters $a_k+r_k$, and final parameter $a_++r_+$, proving (1) for a
maximal reference scale.\\

For an arbitrary positive reference scale $b_n$, use independent
$G_i\sim\Ga(a_i,1)$ and $W_i=b_iG_i/\sum_jb_jG_j$.
The identity $y^{-r_+}=\Gamma(r_+)^{-1}\int_0^\infty
s^{r_+-1}e^{-sy}\dd s$ and Tonelli's theorem give
\[
 \E\!\left[\prod_iW_i^{r_i}\right]
 =\frac{\prod_i b_i^{r_i}(a_i)_{r_i}}{\Gamma(r_+)}
 \int_0^\infty s^{r_+-1}\prod_i(1+b_is)^{-a_i-r_i}\dd s.
\]
Set $s=t/[b_n(1-t)]$. The resulting integral is
\[
 \frac{\prod_{i=1}^{n-1}(b_i/b_n)^{r_i}\prod_i(a_i)_{r_i}}
 {\Gamma(r_+)}
 \int_0^1 t^{r_+-1}(1-t)^{a_+-1}
 \prod_{k=1}^{n-1}(1-z_kt)^{-a_k-r_k}\dd t.
\]
Since $z_k<1$, the integral is finite, and the Euler representation in
Definition~\ref{def:Lauricella} yields (1) with the stated continuation
convention.
\end{proof}

\begin{corollary}\label{cor:carlson}
Let $r_+=\sum_i r_i$ and $a_+=\sum_i a_i$. For shape and scale vectors
$\bm a$ and $\bm b$, Carlson's multivariate $R$-function from
Lemma~\ref{lem:carlson} satisfies
\begin{equation}\tag{2}\label{eq:Carlson_relation}
 R_{-r_+}(\bm a+\bm r;\bm b)
 =b_n^{-r_+}\FD^{(n-1)}\!\left(r_+;
 a_1+r_1,\ldots,a_{n-1}+r_{n-1};a_++r_+;
 1-\frac{b_1}{b_n},\ldots,1-\frac{b_{n-1}}{b_n}\right).
\end{equation}
Hence,
\[
 \boxed{\E\!\left[\prod_{i=1}^nW_i^{r_i}\right]
 =\left(\prod_{i=1}^n b_i^{r_i}\right)
 \frac{\prod_{i=1}^n(a_i)_{r_i}}{(a_+)_{r_+}}
 R_{-r_+}(\bm a+\bm r;\bm b).}
\]
\end{corollary}
\begin{proof}
Apply Lemma~\ref{lem:carlson} under the substitutions
$a\mapsto r_+$, $(b_1,\ldots,b_n)\mapsto(a_1+r_1,\ldots,a_n+r_n)$,
and $(z_1,\ldots,z_n)\mapsto(b_1,\ldots,b_n)$. The corresponding
Lauricella parameter is $c=\sum_i(a_i+r_i)=a_++r_+$. This gives (2),
and substitution into (1) gives the moment representation.
\end{proof}

\begin{corollary}[Mean]\label{cor:mean}
Set $r_i=1$ and $r_j=0$ for all $j\ne i$ in (1). Then
\begin{equation}\tag{3}\label{eq:mean_general}
 \boxed{\begin{aligned}
 \E[W_i]&=\frac{b_i}{b_n}\frac{a_i}{a_+}
 \FD^{(n-1)}\!\bigl(1;
 a_1+\ind_{\{i=1\}},\ldots,a_{n-1}+\ind_{\{i=n-1\}};\\[-2pt]
 &\hspace{43mm}a_++1;z_1,\ldots,z_{n-1}\bigr),\qquad i=1,\ldots,n.
 \end{aligned}}
\end{equation}
\end{corollary}

\begin{corollary}[Second Moment]\label{cor:second}
Set $r_i=2$ and $r_j=0$ for all $j\ne i$ in (1). Then
\begin{equation}\tag{4}\label{eq:second_general}
 \boxed{\begin{aligned}
 \E[W_i^2]&=\left(\frac{b_i}{b_n}\right)^2
 \frac{(a_i)_2}{(a_+)_2}\FD^{(n-1)}\!\bigl(2;
 a_1+2\ind_{\{i=1\}},\ldots,a_{n-1}+2\ind_{\{i=n-1\}};\\[-2pt]
 &\hspace{43mm}a_++2;z_1,\ldots,z_{n-1}\bigr),\qquad i=1,\ldots,n.
 \end{aligned}}
\end{equation}
\end{corollary}

\begin{corollary}[Cross-Moment]\label{cor:cross}
Let $i\ne j$. Set $r_i=r_j=1$ and $r_k=0$ for $k\ne i,j$ in (1). Then
\begin{equation}\tag{5}\label{eq:cross_general}
 \boxed{\begin{aligned}
 \E[W_iW_j]&=\frac{b_i}{b_n}\frac{b_j}{b_n}
 \frac{a_ia_j}{(a_+)_2}\FD^{(n-1)}\!\bigl(2;
 a_1+\ind_{\{i=1\}}+\ind_{\{j=1\}},\ldots,\\[-2pt]
 &\hspace{14mm}a_{n-1}+\ind_{\{i=n-1\}}+\ind_{\{j=n-1\}};
 a_++2;z_1,\ldots,z_{n-1}\bigr).
 \end{aligned}}
\end{equation}
\end{corollary}

When $b_1=\cdots=b_n$, the distribution reduces to the ordinary
Dirichlet, for which $\Cov(W_i,W_j)<0$ for $i\ne j$.
For unequal scales, positive covariances are possible.
For example, numerical evaluation of Equations \eqref{eq:mean_general}, \eqref{eq:second_general} and \eqref{eq:cross_general}, with
$\boldsymbol{a}=(10,10,1)$ and $\boldsymbol{b}=(1,1,100)$ gives $\Cov(W_1,W_2)\approx 0.01125138>0$ with $\Cor(W_1,W_2)\approx 0.7724.$ Positive-correlation examples are also given by \cite{Furman2019}.

\section{Identifiability}
\begin{proposition}[Identifiability of the Scale Parameters]\label{prop:identifiability}
Consider the scaled Dirichlet construction
\[
 W_i=\frac{b_iX_i}{\sum_{j=1}^n b_jX_j},\qquad
 \bm X\sim\Dir_n(\bm a),\quad n\geq2.
\]
The scale parameters are identifiable only up to a common multiplicative
constant: for any $c>0$, $\bm b$ and $\bm b^*=c\bm b$ generate the same
distribution. Conversely, equality of two such distributions identifies
the shape vector and implies proportional scale vectors.
\end{proposition}
\begin{proof}
For $\bm b^*=c\bm b$,
\[
 W_i^*=\frac{b_i^*X_i}{\sum_jb_j^*X_j}
 =\frac{cb_iX_i}{c\sum_jb_jX_j}=W_i.
\]
Thus $\bm W^*=\bm W$ for every realisation of $\bm X$, so common scale
rescaling leaves the distribution unchanged.\\

Conversely, suppose two parameter sets give the same distribution.
As $w_i$ approaches zero with the other coordinates approaching
positive values, the density is proportional to $w_i^{a_i-1}$
to leading order. Equal densities therefore have the same exponent
at each face, so $a_i=\widetilde a_i$ for every $i$.\\

With the shape parameters equal, the density formula gives
\[
 \sum_{i=1}^n \frac{w_i}{b_i}
 = K\sum_{i=1}^n \frac{w_i}{\widetilde b_i}
\]
for a constant $K>0$. Taking limits at each simplex vertex gives
$1/b_i=K/\widetilde b_i$, so $\widetilde b_i=Kb_i$ for every $i$.

\end{proof}

\begin{corollary}[Relative Scale Interpretation]
Multiplying all scale parameters by a common positive constant leaves the
distribution unchanged. The distribution therefore depends on the relative
scales
\[
 \frac{b_1}{b_n},\frac{b_2}{b_n},\ldots,\frac{b_{n-1}}{b_n}.
\]
The scale vector contains one redundant degree of freedom. Setting the
reference scale to one gives the equivalent parameterisation
\[
 \left(a_1,\ldots,a_n,\frac{b_1}{b_n},\ldots,\frac{b_{n-1}}{b_n}\right),
\]
which contains all identifiable information about the distribution.
Consequently, there are $n+(n-1)=2n-1$ identifiable parameters.
\end{corollary}

\section{The Bivariate Scaled Dirichlet Distribution}
Although the general scaled Dirichlet distribution has $n$ components, the case $n=2$ is particularly useful from a computational perspective. In the general case, the mixed moments involve the Lauricella hypergeometric function $F_D^{(n-1)}$. When $n=2$, this reduces to the ordinary Gauss hypergeometric function since
\[
 F_D^{(1)}(\alpha;\beta;\gamma;z)={}_2F_1(\alpha,\beta;\gamma;z).
\]
Let $\mathbf W=(W_1,W_2)\sim\SD_2(\mathbf a,\mathbf b)$, where $\mathbf a=(a_1,a_2)$, $\mathbf b=(b_1,b_2)$, $a_+=a_1+a_2$, and define $z=1-b_1/b_2$. Then
\[
 W_1=\frac{b_1X_1}{b_1X_1+b_2X_2},\qquad
 W_2=\frac{b_2X_2}{b_1X_1+b_2X_2},
 \qquad (X_1,X_2)\sim\Dir_2(a_1,a_2).
\]
Since $W_1+W_2=1$, it follows that $W_2=1-W_1$. Consequently, only one component is free and the density domain of the bivariate scaled Dirichlet distribution is one-dimensional. Therefore, the distribution is completely characterized by the marginal density of $W=W_1$, which is given by
\[
 f_{W_1}(w)=\frac{b_2^{a_1}b_1^{a_2}}{B(a_1,a_2)}
 \frac{w^{a_1-1}(1-w)^{a_2-1}}
 {[b_1(1-w)+b_2w]^{a_1+a_2}},\qquad 0<w<1.
\]
Since $W_2=1-W_1$, substituting $1-w$ into the density of $W_1$ gives
\[
 f_{W_2}(w)=\frac{b_2^{a_1}b_1^{a_2}}{B(a_1,a_2)}
 \frac{w^{a_2-1}(1-w)^{a_1-1}}
 {[b_1w+b_2(1-w)]^{a_1+a_2}},\qquad 0<w<1.
\]
The hypergeometric functions below are understood through the Euler-integral continuation when $|z|\geq1$, as specified in Theorem~1 and Definition~2.

\subsection{Moments}
\subsubsection{First Moments}
Using Equation~\eqref{eq:mean_general}, the first moments are
\[
 \E[W_1]=\frac{b_1}{b_2}\frac{a_1}{a_+}
 {}_2F_1(1,a_1+1;a_++1;z)
\]
and
\[
 \E[W_2]=\frac{a_2}{a_+}{}_2F_1(1,a_1;a_++1;z).
\]

\subsubsection{Second Moments}
Using Equation~\eqref{eq:second_general}, the second moments are
\[
 \E[W_1^2]=\left(\frac{b_1}{b_2}\right)^2\frac{(a_1)_2}{(a_+)_2}
 {}_2F_1(2,a_1+2;a_++2;z)
\]
and
\[
 \E[W_2^2]=\frac{(a_2)_2}{(a_+)_2}{}_2F_1(2,a_1;a_++2;z).
\]

\subsubsection{Variances}
Utilising the first and second moments, we may obtain
$\Var(W_1)=\E[W_1^2]-(\E[W_1])^2$ and
$\Var(W_2)=\E[W_2^2]-(\E[W_2])^2$.
Additionally, since $W_1+W_2=1$, it also follows that $\Var(W_1)=\Var(W_2)$.

\subsubsection{Covariance \& Correlation}
Utilising Equation~\eqref{eq:cross_general}, the cross-moment is
\[
 \E[W_1W_2]=\frac{b_1}{b_2}\frac{a_1a_2}{(a_+)_2}
 {}_2F_1(2,a_1+1;a_++2;z).
\]
Hence, we may obtain $\Cov(W_1,W_2)=\E[W_1W_2]-\E[W_1]\E[W_2]$.
Furthermore, since $W_2=1-W_1$, this covariance can also be written as
$\Cov(W_1,W_2)=-\Var(W_1)=-\Var(W_2)$. Additionally,
$\Cor(W_1,W_2)=-1$, which implies that the two components are perfectly negatively correlated and, therefore, one component is completely determined by the other. Unequal scales nevertheless extend the marginal family beyond the ordinary Beta distribution. The bivariate case therefore adds marginal flexibility, but cannot change this deterministic correlation constraint.

\subsubsection{A Dependent-Gamma Extension}

The scaled Dirichlet distribution admits an equivalent representation in terms of normalized independent Gamma random variables. Let $Y_i \stackrel{\mathrm{ind}}{\sim}
\mathrm{Gamma}\left(a_i,\frac{1}{b_i}\right)$ for $i=1,\ldots,n$
where the second parameter denotes the rate. Then $W_i=\frac{Y_i}{\sum_{j=1}^nY_j}$ for $ i=1,\ldots,n$
implies that $\mathbf{W}=(W_1,\ldots,W_n)^\top$ follows a scaled Dirichlet distribution.\\

This representation suggests a natural extension obtained by relaxing the assumption of independence among the underlying Gamma random variables. In the bivariate case, let $Y_1$ have a Gamma marginal distribution and specify the conditional distribution of $Y_2$ given $Y_1>y_1$ using the accelerated failure conditionals framework of \citet{lakhani2025modelsacceleratedfailureconditionals}. For example, suppose $Y_1\sim \mathrm{Gamma}\left(a_1,\frac{1}{b_1}\right)$ and $Y_2\mid Y_1>y_1
\sim
\mathrm{Gamma}\left(a_2,\beta(y_1)\right)$,
where $\beta(y_1)$ denotes the acceleration function responsible for ensuring a dependence relationship between the two variables. We assume that $\beta$ is chosen so that these conditional
distributions and the marginal distribution of $Y_1$ define
a valid joint density. Furthermore, define $W_1=\frac{Y_1}{Y_1+Y_2}$ and $W_2=\frac{Y_2}{Y_1+Y_2}$. Now, the resulting random vector, $\mathbf{W}$, remains supported on the simplex since $W_1+W_2=1$, but it does not generally follow a scaled Dirichlet distribution because the Gamma constituents are no longer independent. Instead, this construction defines a dependent-Gamma extension of the scaled Dirichlet family.\\

To derive the distribution of $W_1$, introduce the transformation $W=\frac{Y_1}{Y_1+Y_2}$ with $T=Y_1+Y_2$ and inverse transformation $Y_1=WT$ and $Y_2=(1-W)T$. The absolute value of the Jacobian determinant is $\left|\frac{\partial(Y_1,Y_2)}{\partial(W,T)}\right|=T$.
Therefore $f_{W,T}(w,t)
= f_{Y_1, Y_2}\left(wt, (1-w)t \right)t$
and the marginal density of $W$ is $f_W(w) = 
\int_0^\infty f_{Y_1, Y_2}\left(wt, (1-w)t \right) dt$ with $0<w<1$.\\

When the acceleration function is constant,\ $\beta(y_1)=\frac{1}{b_2}$, the conditional distribution of $Y_2$ no longer depends on $Y_1$, so that $Y_1$ and $Y_2$ are independent Gamma random variables. The standard bivariate scaled Dirichlet distribution is therefore recovered as a special case of the accelerated failure conditionals model given in \citet{lakhani2025modelsacceleratedfailureconditionals}.

\section{The Trivariate Scaled Dirichlet Distribution}
Although the general scaled Dirichlet distribution has $n$ components, the case $n=3$ is particularly useful from both a computational and modelling perspective. In the general case, the mixed moments involve the Lauricella hypergeometric function $F_D^{(n-1)}$. When $n=3$, this reduces to Appell's hypergeometric function $F_1$ since
\[
 F_D^{(2)}(\alpha;\beta_1,\beta_2;\gamma;z_1,z_2)
 =F_1(\alpha;\beta_1,\beta_2;\gamma;z_1,z_2).
\]
Let $\mathbf W=(W_1,W_2,W_3)\sim\SD_3(\mathbf a,\mathbf b)$, where $\mathbf a=(a_1,a_2,a_3)$, $\mathbf b=(b_1,b_2,b_3)$, $a_+=a_1+a_2+a_3$, and define $z_1=1-b_1/b_3$, $z_2=1-b_2/b_3$. Then
\[
 W_1=\frac{b_1X_1}{b_1X_1+b_2X_2+b_3X_3},\quad
 W_2=\frac{b_2X_2}{b_1X_1+b_2X_2+b_3X_3},\quad
 W_3=\frac{b_3X_3}{b_1X_1+b_2X_2+b_3X_3},
\]
where $(X_1,X_2,X_3)\sim\Dir_3(a_1,a_2,a_3)$. Since $W_1+W_2+W_3=1$, it is sufficient to write the joint density of $(W_1,W_2)$, with $W_3=1-W_1-W_2$. The density is
\[
 f_{W_1,W_2}(w_1,w_2)=
 \frac{b_1^{a_+-a_1}b_2^{a_+-a_2}b_3^{a_+-a_3}}{B(a_1,a_2,a_3)}
 \frac{w_1^{a_1-1}w_2^{a_2-1}(1-w_1-w_2)^{a_3-1}}
 {[b_2b_3w_1+b_1b_3w_2+b_1b_2(1-w_1-w_2)]^{a_+}},
\]
where $w_1>0$, $w_2>0$, and $w_1+w_2<1$.\\

Now, we define the relative scale parameters $\rho_1=b_1/b_3$ and $\rho_2=b_2/b_3$. The density is then
\begin{equation}\tag{6}\label{eq:tri_density}
 f_{W_1,W_2}(w_1,w_2)=
 \frac{1}{B(a_1,a_2,a_3)\rho_1^{a_1}\rho_2^{a_2}}
 \frac{w_1^{a_1-1}w_2^{a_2-1}(1-w_1-w_2)^{a_3-1}}
 {[w_1/\rho_1+w_2/\rho_2+(1-w_1-w_2)]^{a_+}},
\end{equation}
where $w_1>0$, $w_2>0$, and $w_1+w_2<1$.

\begin{proposition}[Marginal Density of $W_3$]\label{prop:marginal_W3}
Let $(W_1,W_2,W_3)\sim\SD_3(\mathbf a,\mathbf b)$, and define the relative scale parameters $\rho_1=b_1/b_3$ and $\rho_2=b_2/b_3$. Let $a_+=a_1+a_2+a_3$. Then the marginal density of $W_3$ is
\[
 f_{W_3}(t)=\frac{t^{a_3-1}(1-t)^{a_1+a_2-1}}
 {B(a_1+a_2,a_3)\rho_1^{a_1}\rho_2^{a_2}}
 A(t)^{-a_+}{}_2F_1(a_+,a_1;a_1+a_2;x(t)),
\]
where $0<t<1$,
\[
 A(t)=t+\frac{1-t}{\rho_2},\qquad
 x(t)=\frac{(1-t)(\rho_2^{-1}-\rho_1^{-1})}{t+(1-t)/\rho_2}.
\]
\end{proposition}
\begin{proof}
The joint density of $(W_1,W_2)$ is given by Equation~\eqref{eq:tri_density}. To obtain the marginal density of $W_3$, set $t=w_3$. Since $w_1+w_2+w_3=1$, fixing $w_3=t$ implies $w_1+w_2=1-t$. We parameterise this split by introducing $u=w_1/(1-t)$. Hence $w_1=(1-t)u$ and $w_2=(1-t)(1-u)$. The conditions $w_1>0$ and $w_2>0$ imply $0<u<1$. The absolute Jacobian for $(u,t)\mapsto(w_1,w_2)$ is $1-t$. Therefore,
\[
 f_{W_3}(t)=\int_0^1 f_{W_1,W_2}((1-t)u,(1-t)(1-u))(1-t)\,\mathrm du.
\]
Substituting $w_1=(1-t)u$ and $w_2=(1-t)(1-u)$, we obtain
\[
 f_{W_3}(t)=\frac{t^{a_3-1}(1-t)^{a_1+a_2-1}}
 {B(a_1,a_2,a_3)\rho_1^{a_1}\rho_2^{a_2}}
 \int_0^1\frac{u^{a_1-1}(1-u)^{a_2-1}}
 {[t+(1-t)(u/\rho_1+(1-u)/\rho_2)]^{a_+}}\,\mathrm du.
\]
Now write the denominator inside the integral as
\[
 t+(1-t)\left(\frac{u}{\rho_1}+\frac{1-u}{\rho_2}\right)
 =A(t)(1-x(t)u).
\]
Thus,
\[
 \left[t+(1-t)\left(\frac{u}{\rho_1}+\frac{1-u}{\rho_2}\right)\right]^{-a_+}
 =A(t)^{-a_+}(1-x(t)u)^{-a_+}.
\]
Hence,
\[
 f_{W_3}(t)=\frac{t^{a_3-1}(1-t)^{a_1+a_2-1}}
 {B(a_1,a_2,a_3)\rho_1^{a_1}\rho_2^{a_2}}A(t)^{-a_+}
 \int_0^1 u^{a_1-1}(1-u)^{a_2-1}(1-x(t)u)^{-a_+}\,\mathrm du.
\]
Using the Euler integral representation of the Gauss hypergeometric function,
\[
 \int_0^1u^{a_1-1}(1-u)^{a_2-1}(1-xu)^{-a_+}\,\mathrm du
 =B(a_1,a_2){}_2F_1(a_+,a_1;a_1+a_2;x),
\]
we get
\[
 f_{W_3}(t)=\frac{B(a_1,a_2)t^{a_3-1}(1-t)^{a_1+a_2-1}}
 {B(a_1,a_2,a_3)\rho_1^{a_1}\rho_2^{a_2}}A(t)^{-a_+}
 {}_2F_1(a_+,a_1;a_1+a_2;x(t)).
\]
Finally, using $B(a_1,a_2)/B(a_1,a_2,a_3)=1/B(a_1+a_2,a_3)$, we obtain the stated density. The Euler representation is valid here because $A(t)>0$ and
\[
 1-x(t)=\frac{t+(1-t)/\rho_1}{A(t)}>0.
\]
\end{proof}

\subsection{Moments}
The Appell $F_1$ series converges when $|z_1|,|z_2|<1$.
Positive scales ensure that $z_1,z_2<1$, but either argument can
be less than $-1$. In that case, $F_1$ can be as specified in Theorem \ref{thm:mixed} and Definition \ref{def:Lauricella}.

\subsubsection{First Moments}
Using Equation~\eqref{eq:mean_general}, the first moments are
\begin{equation}\tag{7}\label{eq:tri_mean1}
 \E[W_1]=\rho_1\frac{a_1}{a_+}F_1(1;a_1+1,a_2;a_++1;z_1,z_2),
\end{equation}
\begin{equation}\tag{8}\label{eq:tri_mean2}
 \E[W_2]=\rho_2\frac{a_2}{a_+}F_1(1;a_1,a_2+1;a_++1;z_1,z_2),
\end{equation}
and
\[
 \E[W_3]=\frac{a_3}{a_+}F_1(1;a_1,a_2;a_++1;z_1,z_2).
\]

\subsubsection{Second Moments}
Using Equation~\eqref{eq:second_general}, the second moments are
\begin{equation}\tag{9}\label{eq:tri_second1}
 \E[W_1^2]=\rho_1^2\frac{(a_1)_2}{(a_+)_2}
 F_1(2;a_1+2,a_2;a_++2;z_1,z_2),
\end{equation}
\begin{equation}\tag{10}\label{eq:tri_second2}
 \E[W_2^2]=\rho_2^2\frac{(a_2)_2}{(a_+)_2}
 F_1(2;a_1,a_2+2;a_++2;z_1,z_2),
\end{equation}
and
\[
 \E[W_3^2]=\frac{(a_3)_2}{(a_+)_2}F_1(2;a_1,a_2;a_++2;z_1,z_2).
\]

\subsubsection{Variances}
Utilising the first and second moments, the variances may be computed as
$\Var(W_1)=\E[W_1^2]-(\E[W_1])^2$,
$\Var(W_2)=\E[W_2^2]-(\E[W_2])^2$, and
$\Var(W_3)=\E[W_3^2]-(\E[W_3])^2$.

\subsubsection{Covariances \& Correlations}
Using Equation~\eqref{eq:cross_general}, the cross-moments are
\begin{equation}\tag{11}\label{eq:tri_cross12}
 \E[W_1W_2]=\rho_1\rho_2\frac{a_1a_2}{(a_+)_2}
 F_1(2;a_1+1,a_2+1;a_++2;z_1,z_2),
\end{equation}
\[
 \E[W_1W_3]=\rho_1\frac{a_1a_3}{(a_+)_2}
 F_1(2;a_1+1,a_2;a_++2;z_1,z_2),
\]
and
\[
 \E[W_2W_3]=\rho_2\frac{a_2a_3}{(a_+)_2}
 F_1(2;a_1,a_2+1;a_++2;z_1,z_2).
\]
Hence, we may compute $\Cov(W_i,W_j)=\E[W_iW_j]-\E[W_i]\E[W_j]$, for $i\ne j$. Unlike the bivariate case, the simplex constraint $W_1+W_2+W_3=1$ does not force any single pair of components to be perfectly negatively correlated. Instead, for positive shapes and scales, the pairwise correlations satisfy $-1<\Cor(W_i,W_j)<1$ for $i\ne j$.

\subsection{Method-of-Moments Estimators}
For $n=3$, the identifiable parameter vector is $\boldsymbol\theta=(a_1,a_2,a_3,\rho_1,\rho_2)$. Given observations $\mathbf w_1,\ldots,\mathbf w_N$ with $\mathbf w_\ell=(w_{1\ell},w_{2\ell},w_{3\ell})$, define the empirical moments
\[
 \widehat m_1=\frac1N\sum_{\ell=1}^Nw_{1\ell},\qquad
 \widehat m_2=\frac1N\sum_{\ell=1}^Nw_{2\ell},
\]
\[
 \widehat m_{11}=\frac1N\sum_{\ell=1}^Nw_{1\ell}^2,\qquad
 \widehat m_{22}=\frac1N\sum_{\ell=1}^Nw_{2\ell}^2,
 \qquad
 \widehat m_{12}=\frac1N\sum_{\ell=1}^Nw_{1\ell}w_{2\ell}.
\]
The method-of-moments estimator is obtained, when an admissible solution exists, utilising Equations~\eqref{eq:tri_mean1}, \eqref{eq:tri_mean2}, \eqref{eq:tri_second1}, \eqref{eq:tri_second2}, and~\eqref{eq:tri_cross12} and numerically solving
\[
 \widehat m_1=\E[W_1],\quad \widehat m_2=\E[W_2],\quad
 \widehat m_{11}=\E[W_1^2],\quad \widehat m_{22}=\E[W_2^2],\quad
 \widehat m_{12}=\E[W_1W_2].
\]

\subsection{Maximum Likelihood Estimators}
Let $\boldsymbol\theta=(a_1,a_2,a_3,\rho_1,\rho_2)$ and $a_+=a_1+a_2+a_3$. For an independent random sample $(w_{11},w_{21}),\ldots,(w_{1N},w_{2N})$, the likelihood function is
\[
 L(\boldsymbol\theta)=\prod_{i=1}^Nf_{W_1,W_2}(w_{1i},w_{2i};\boldsymbol\theta).
\]
Therefore, using Equation~\eqref{eq:tri_density}, we have
\[
 L(\boldsymbol\theta)=\prod_{i=1}^N
 \frac{w_{1i}^{a_1-1}w_{2i}^{a_2-1}(1-w_{1i}-w_{2i})^{a_3-1}}
 {B(a_1,a_2,a_3)\rho_1^{a_1}\rho_2^{a_2}
 [w_{1i}/\rho_1+w_{2i}/\rho_2+(1-w_{1i}-w_{2i})]^{a_+}}.
\]
Taking logs, the log-likelihood is
\begin{align*}
 \ell(\boldsymbol\theta)={}&-N\log B(a_1,a_2,a_3)-Na_1\log\rho_1-Na_2\log\rho_2
 +(a_1-1)\sum_{i=1}^N\log w_{1i}\\
 &+(a_2-1)\sum_{i=1}^N\log w_{2i}
 +(a_3-1)\sum_{i=1}^N\log(1-w_{1i}-w_{2i})\\
 &-a_+\sum_{i=1}^N\log\left[\frac{w_{1i}}{\rho_1}+\frac{w_{2i}}{\rho_2}+(1-w_{1i}-w_{2i})\right].
\end{align*}
The maximum likelihood estimators of $\boldsymbol\theta=(a_1,a_2,a_3,\rho_1,\rho_2)$ are obtained by maximizing the log-likelihood function,
$\widehat{\boldsymbol\theta}=\operatorname*{arg\,max}_{\boldsymbol\theta}\ell(\boldsymbol\theta)$, subject to $a_1,a_2,a_3,\rho_1,\rho_2>0$.

\subsection{Likelihood Ratio Test Against the Ordinary Dirichlet}
The general form of the likelihood ratio test statistic is given by
\[
 \Lambda=\frac{\sup_{\boldsymbol\theta\in\Theta_0}L(\boldsymbol\theta)}
 {\sup_{\boldsymbol\theta\in\Theta}L(\boldsymbol\theta)},
\]
where $\Theta$ denotes the unrestricted parameter space and $\Theta_0\subset\Theta$ denotes the parameter space under the null hypothesis. For independent and identically distributed observations, under the null hypothesis and the regularity conditions for Wilks' theorem, $-2\log\Lambda\xrightarrow{d}\chi^2_{\dim(\Theta)-\dim(\Theta_0)}$. Thus, under these assumptions, the null hypothesis is rejected at significance level $\alpha$ whenever
\[
 -2\log\Lambda>\chi^2_{\dim(\Theta)-\dim(\Theta_0),1-\alpha},
\]
where the right-hand side denotes the $(1-\alpha)$-quantile.\\

The ordinary Dirichlet distribution is obtained as a special case of the scaled Dirichlet distribution under the restriction $\rho_1=\rho_2=1$. Under this restriction, the density from Equation~\eqref{eq:tri_density} becomes
\[
 f_D(w_1,w_2;\mathbf a)=\frac{w_1^{a_1-1}w_2^{a_2-1}(1-w_1-w_2)^{a_3-1}}
 {B(a_1,a_2,a_3)}.
\]
The likelihood ratio test considers
\[
 H_0:\rho_1=\rho_2=1
 \qquad\text{against}\qquad
 H_1:\text{at least one of }\rho_1,\rho_2\ne1.
\]
Rejection of the null hypothesis provides evidence that allowing unequal scale parameters improves the fit relative to the ordinary Dirichlet distribution. Conversely, failure to reject provides insufficient evidence against the equal-scale restriction; it does not establish that the ordinary Dirichlet is an adequate model.\\

Let $\Theta=\{(a_1,a_2,a_3,\rho_1,\rho_2)^{\mathsf T}:a_1,a_2,a_3,\rho_1,\rho_2>0\}$, with $(\widehat{\mathbf a}_S,\widehat{\boldsymbol\rho}_S)$ denoting the unrestricted maximum likelihood estimator. Let $\Theta_0=\{(a_1,a_2,a_3,1,1)^{\mathsf T}:a_1,a_2,a_3>0\}$, with $\widehat{\mathbf a}_D$ denoting the restricted maximum likelihood estimator. The likelihood ratio is
\[
 \Lambda=\frac{L_D(\widehat{\mathbf a}_D)}
 {L_S(\widehat{\mathbf a}_S,\widehat{\boldsymbol\rho}_S)}
 =\frac{\displaystyle\prod_{i=1}^N
 \dfrac{w_{1i}^{\widehat a_{1,D}-1}w_{2i}^{\widehat a_{2,D}-1}
 (1-w_{1i}-w_{2i})^{\widehat a_{3,D}-1}}{B(\widehat{\mathbf a}_D)}}
 {\displaystyle\prod_{i=1}^N
 \dfrac{w_{1i}^{\widehat a_{1,S}-1}w_{2i}^{\widehat a_{2,S}-1}
 (1-w_{1i}-w_{2i})^{\widehat a_{3,S}-1}}
 {B(\widehat{\mathbf a}_S)\widehat\rho_{1,S}^{\widehat a_{1,S}}
 \widehat\rho_{2,S}^{\widehat a_{2,S}}
 [w_{1i}/\widehat\rho_{1,S}+w_{2i}/\widehat\rho_{2,S}+(1-w_{1i}-w_{2i})]^{\widehat a_{+,S}}}}.
\]
The likelihood ratio statistic is then
\[
 -2\log\Lambda=2[\ell_S(\widehat{\mathbf a}_S,\widehat{\boldsymbol\rho}_S)
 -\ell_D(\widehat{\mathbf a}_D)],
\]
where $\ell_S=\log L_S$ and $\ell_D=\log L_D$. Under $H_0$ and the assumptions stated above, $-2\log\Lambda\xrightarrow{d}\chi_2^2$, with $\dim(\Theta)-\dim(\Theta_0)=2$, since the scaled Dirichlet model contains two additional identifiable parameters: $\rho_1$ and $\rho_2$.

\subsection{Simulations}
Simulating observations $\mathbf w=(w_1,w_2,w_3)\sim\SD_3(\mathbf a,\mathbf b)$ is fairly straightforward and is given in Algorithm~1. For each sample size, we generated $R=\mathrm{10, 000}$ simulated
samples using Algorithm~1.
\begin{center}
\begin{minipage}{0.96\linewidth}
\hrule\smallskip
\textbf{Algorithm 1}\quad Simulation from $\SD_3(\mathbf a,\mathbf b)$
\smallskip\hrule
\begin{enumerate}
 \item Generate independent random variables
 $G_i\sim\Ga(a_i,1)$, $i=1,2,3$.
 \item Set $W_i=b_iG_i/\sum_{j=1}^3b_jG_j$, $i=1,2,3$.
 \item Then $\mathbf W=(W_1,W_2,W_3)\sim\SD_3(\mathbf a,\mathbf b)$.
\end{enumerate}
\hrule
\end{minipage}
\end{center}

\begin{table}[htbp]
\centering
\resizebox{\textwidth}{!}{%
\begin{tabular}{clrrrrrrr}
\toprule
$N$ & Parameter & MME & MLE & SE(MME) & SE(MLE) & $\mathrm{PC}(W_1,W_2)$ & $\mathrm{PC}(W_1,W_3)$ & $\mathrm{PC}(W_2,W_3)$\\
\midrule
\multirow{8}{*}{500} & $a_1$ & $0.099$ & $0.100$ & $0.057$ & $0.005$ &  &  &  \\
 & $a_2$ & $0.199$ & $0.201$ & $0.298$ & $0.012$ &  &  &  \\
 & $a_3$ & $0.300$ & $0.301$ & $0.141$ & $0.023$ &  &  &  \\
 & $\rho_1$ & $0.408$ & $0.403$ & $4.135$ & $0.159$ &  &  &  \\
 & $\rho_2$ & $0.510$ & $0.502$ & $3.712$ & $0.189$ &  &  &  \\
 & $\Cor(W_1,W_2)$ & $-0.238$ & $-0.238$ & $0.023$ & $0.019$ & $-0.238$ &  &  \\
 & $\Cor(W_1,W_3)$ & $-0.458$ & $-0.458$ & $0.030$ & $0.025$ &  & $-0.458$ &  \\
 & $\Cor(W_2,W_3)$ & $-0.754$ & $-0.754$ & $0.026$ & $0.023$ &  &  & $-0.754$ \\
\midrule
\multirow{8}{*}{1000} & $a_1$ & $0.100$ & $0.100$ & $0.019$ & $0.004$ &  &  &  \\
 & $a_2$ & $0.200$ & $0.200$ & $0.043$ & $0.009$ &  &  &  \\
 & $a_3$ & $0.301$ & $0.301$ & $0.033$ & $0.016$ &  &  &  \\
 & $\rho_1$ & $0.404$ & $0.401$ & $0.316$ & $0.108$ &  &  &  \\
 & $\rho_2$ & $0.509$ & $0.502$ & $0.410$ & $0.128$ &  &  &  \\
 & $\Cor(W_1,W_2)$ & $-0.238$ & $-0.238$ & $0.016$ & $0.013$ & $-0.238$ &  &  \\
 & $\Cor(W_1,W_3)$ & $-0.459$ & $-0.459$ & $0.021$ & $0.018$ &  & $-0.459$ &  \\
 & $\Cor(W_2,W_3)$ & $-0.754$ & $-0.753$ & $0.018$ & $0.016$ &  &  & $-0.754$ \\
\midrule
\multirow{8}{*}{2000} & $a_1$ & $0.100$ & $0.100$ & $0.012$ & $0.003$ &  &  &  \\
 & $a_2$ & $0.200$ & $0.200$ & $0.022$ & $0.006$ &  &  &  \\
 & $a_3$ & $0.300$ & $0.300$ & $0.023$ & $0.011$ &  &  &  \\
 & $\rho_1$ & $0.405$ & $0.400$ & $0.192$ & $0.075$ &  &  &  \\
 & $\rho_2$ & $0.500$ & $0.500$ & $0.245$ & $0.088$ &  &  &  \\
 & $\Cor(W_1,W_2)$ & $-0.238$ & $-0.238$ & $0.012$ & $0.009$ & $-0.238$ &  &  \\
 & $\Cor(W_1,W_3)$ & $-0.459$ & $-0.459$ & $0.015$ & $0.013$ &  & $-0.459$ &  \\
 & $\Cor(W_2,W_3)$ & $-0.754$ & $-0.754$ & $0.013$ & $0.011$ &  &  & $-0.754$ \\
\bottomrule
\end{tabular}}
\caption{Simulation results for the trivariate scaled Dirichlet distribution with true parameter values $a_1=0.1$, $a_2=0.2$, $a_3=0.3$, $\rho_1=0.4$, $\rho_2=0.5$. The columns labeled MME and MLE report the method-of-moments and maximum likelihood estimates, respectively. The columns labeled SE(MME) and SE(MLE) report the Monte Carlo sampling standard deviations of these estimators. The PC columns report the empirical Pearson correlations for each pair of variables.}
\label{tab:simulation}
\end{table}

Table~\ref{tab:simulation} shows that both methods give average
estimates close to the true parameter values. However, the maximum
likelihood estimates vary less across simulated samples. This
difference is most noticeable for the relative scales $\rho_1$ and
$\rho_2$ when $N=500$. The Monte Carlo standard deviations\footnote{The reported standard errors are Monte Carlo estimates of the sampling standard deviation of the estimator, computed over $R$ simulation replications. Specifically, for an estimator $\widehat\theta$, $\widehat{\operatorname{SE}}(\widehat\theta)
 =\sqrt{\frac{1}{R-1}\sum_{r=1}^R
 (\widehat\theta^{(r)}-\overline{\widehat\theta})^2}$, where $\widehat\theta^{(r)}$ denotes the estimate from the $r$th replication and $\overline{\widehat\theta}$ is the average estimate across all replications. } decrease
for both methods as the sample size increases. The two methods also
give nearly identical average estimates of the pairwise correlations,
which are close to the corresponding sample Pearson correlations.

\subsection{An Example Dataset}
We illustrate the trivariate scaled Dirichlet distribution using
daily closing prices for the E-mini S\&P 500, E-mini Nasdaq-100, and
E-mini Dow Jones Industrial Average futures. We downloaded the
\texttt{ES=F}, \texttt{NQ=F}, and \texttt{YM=F} series from Yahoo Finance
\citep{yahoofinance2026} using the R package \texttt{quantmod}, with a
requested date range from 1 January 2010 to 2 July 2026. We kept dates
with closing prices for all three series and calculated log returns
between consecutive dates in the resulting sample.\\

For each date $t$, we calculated the share of the sum of squared
returns attributable to each futures series:
\[
 W_{i,t}=\frac{R_{i,t}^2}{\sum_{j=1}^3 R_{j,t}^2},
 \qquad i=1,2,3,
\]
where $R_{i,t}$ is the log return for series $i$. We excluded days
with a zero denominator and retained only observations for which
all three shares were finite and greater than $10^{-8}$. The
final sample contained 4095 observations from 2010-01-05 
to 2026-07-01.\\

Table~\ref{tab:futures} compares the ordinary and scaled Dirichlet
fits. The scaled Dirichlet has a higher maximized log-likelihood
and a lower AIC, even after accounting for its two additional
scale parameters. Its maximum likelihood fit comes closer to the
sample Pearson correlations for the ES-YM and NQ-YM pairs, while
the ordinary Dirichlet is closer for ES--NQ. The likelihood-ratio statistic is $597.162$. We report it as a
descriptive comparison of the two fits because independence
across trading days was not assessed.

\begin{table}[htbp]
\centering
\setlength{\tabcolsep}{3pt}
\resizebox{\textwidth}{!}{%
\begin{tabular}{clrrrrrrrrrrrr}
\toprule
Model & Parameter & MME & MLE & $\widehat m_1$ & $\widehat m_2$ & $\widehat m_{11}$ & $\widehat m_{22}$ & $\widehat m_{12}$ & $\mathrm{PC}(W_1,W_2)$ & $\mathrm{PC}(W_1,W_3)$ & $\mathrm{PC}(W_2,W_3)$ & $\ell$ & AIC\\
\midrule
\multirow{13}{*}{\rotatebox{90}{Scaled Dirichlet}} & $a_1$ & $16.261$ & $1.209$ &  &  &  &  &  &  &  &  & \multirow{13}{*}{$3449.293$} & \multirow{13}{*}{$-6888.585$} \\
 & $a_2$ & $0.802$ & $0.754$ &  &  &  &  &  &  &  &  &  &  \\
 & $a_3$ & $1.212$ & $0.859$ &  &  &  &  &  &  &  &  &  &  \\
 & $\rho_1$ & $0.053$ & $0.647$ &  &  &  &  &  &  &  &  &  &  \\
 & $\rho_2$ & $3.233$ & $2.326$ &  &  &  &  &  &  &  &  &  &  \\
 & $\E[W_1]$ & $0.262$ & $0.280$ & $0.263$ &  &  &  &  &  &  &  &  &  \\
 & $\E[W_2]$ & $0.451$ & $0.441$ &  & $0.448$ &  &  &  &  &  &  &  &  \\
 & $\E[W_1^2]$ & $0.096$ & $0.128$ &  &  & $0.095$ &  &  &  &  &  &  &  \\
 & $\E[W_2^2]$ & $0.271$ & $0.271$ &  &  &  & $0.273$ &  &  &  &  &  &  \\
 & $\E[W_1W_2]$ & $0.092$ & $0.087$ &  &  &  &  & $0.094$ &  &  &  &  &  \\
 & $\Cor(W_1,W_2)$ & $-0.621$ & $-0.590$ &  &  &  &  &  & $-0.537$ &  &  &  &  \\
 & $\Cor(W_1,W_3)$ & $-0.014$ & $-0.257$ &  &  &  &  &  &  & $-0.071$ &  &  &  \\
 & $\Cor(W_2,W_3)$ & $-0.775$ & $-0.629$ &  &  &  &  &  &  &  & $-0.803$ &  &  \\
\midrule
\multirow{11}{*}{\rotatebox{90}{Dirichlet}} & $a_1$ & $0.930$ & $0.810$ &  &  &  &  &  &  &  &  & \multirow{11}{*}{$3150.712$} & \multirow{11}{*}{$-6295.424$} \\
 & $a_2$ & $1.658$ & $1.075$ &  &  &  &  &  &  &  &  &  &  \\
 & $a_3$ & $1.034$ & $0.750$ &  &  &  &  &  &  &  &  &  &  \\
 & $\E[W_1]$ & $0.257$ & $0.307$ & $0.263$ &  &  &  &  &  &  &  &  &  \\
 & $\E[W_2]$ & $0.458$ & $0.407$ &  & $0.448$ &  &  &  &  &  &  &  &  \\
 & $\E[W_1^2]$ & $0.107$ & $0.152$ &  &  & $0.095$ &  &  &  &  &  &  &  \\
 & $\E[W_2^2]$ & $0.263$ & $0.232$ &  &  &  & $0.273$ &  &  &  &  &  &  \\
 & $\E[W_1W_2]$ & $0.092$ & $0.090$ &  &  &  &  & $0.094$ &  &  &  &  &  \\
 & $\Cor(W_1,W_2)$ & $-0.540$ & $-0.550$ &  &  &  &  &  & $-0.537$ &  &  &  &  \\
 & $\Cor(W_1,W_3)$ & $-0.372$ & $-0.422$ &  &  &  &  &  &  & $-0.071$ &  &  &  \\
 & $\Cor(W_2,W_3)$ & $-0.581$ & $-0.525$ &  &  &  &  &  &  &  & $-0.803$ &  &  \\
\bottomrule
\end{tabular}}
\caption{Parameter estimation and model comparison for the futures data under the scaled Dirichlet and Dirichlet models. The MME and MLE columns report the method-of-moments and maximum likelihood estimates, respectively. The columns $\widehat m_1$, $\widehat m_2$, $\widehat m_{11}$, $\widehat m_{22}$, and $\widehat m_{12}$ contain the corresponding empirical moments. The PC columns report the empirical Pearson correlations for each pair of variables. The final two columns give the maximized log-likelihood ($\ell$) and the Akaike Information Criterion (AIC).}
\label{tab:futures}
\end{table}

The figures below illustrate the fitted joint and marginal densities using the maximum likelihood estimates in Table~\ref{tab:futures}.

\clearpage
\begin{figure}[p]
 \centering
 \includegraphics[width=\textwidth]{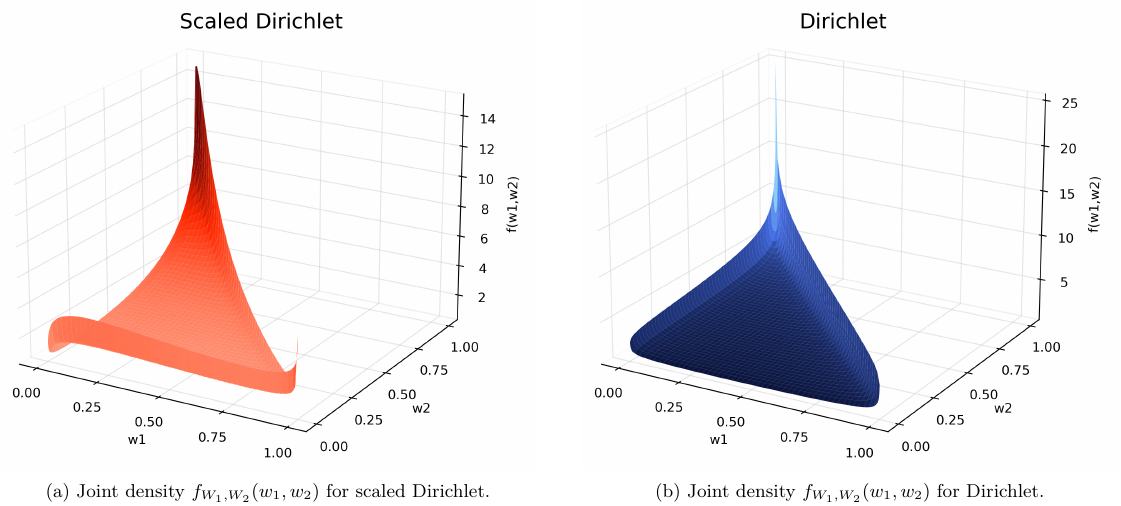}
 \caption{Estimated bivariate joint density, $f_{W_1,W_2}(w_1,w_2)$, under the scaled Dirichlet and Dirichlet models fitted to the futures data using maximum likelihood estimators from Table~\ref{tab:futures}.}
 \label{fig:joint_density}
\end{figure}
\begin{figure}[p]
 \centering
 \includegraphics[width=\textwidth]{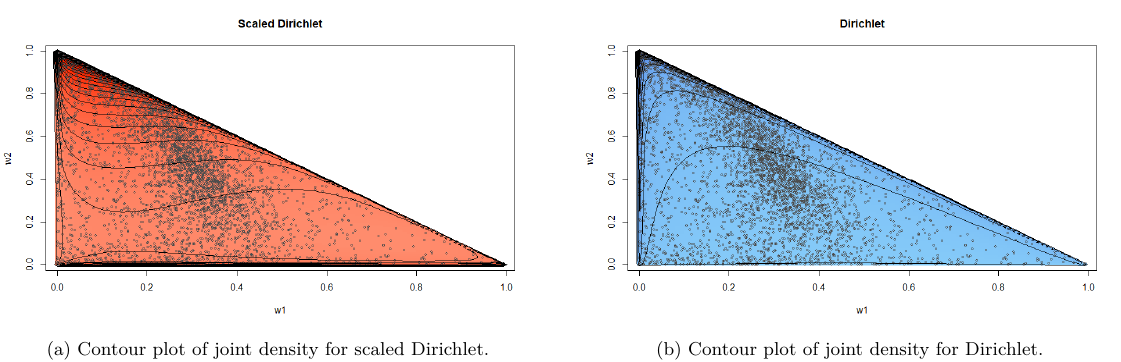}
 \caption{Contour plots of estimated bivariate joint density, $f_{W_1,W_2}(w_1,w_2)$, under the scaled Dirichlet and Dirichlet models fitted to the futures data using maximum likelihood estimators from Table~\ref{tab:futures}. The scatter points represent the observed $(W_1,W_2)$ pairs, while the contour lines correspond to level curves of the fitted joint density.}
 \label{fig:contours}
\end{figure}
\clearpage
\begin{figure}[!ht]
 \centering
 \includegraphics[width=0.48\textwidth]{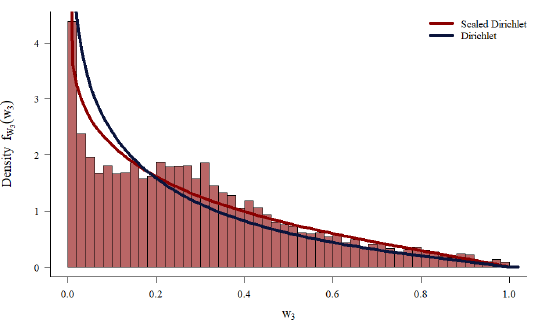}
 \caption{Estimated marginal density, $f_{W_3}(w_3)$, under the scaled Dirichlet and Dirichlet models (red and blue curves, respectively). The histogram shows the empirical distribution of the observed values of $W_3$.}
 \label{fig:marginal_W3}
\end{figure}

\appendix
\section{Additional Lemmas and Definitions}
\begin{lemma}[Dirichlet Mixed Moments; \cite{ng2011dirichlet}]\label{lem:dirichlet_moments}
Let $\bm X=(X_1,\ldots,X_n)^{\mathsf T}\sim\Dir_n(a_1,\ldots,a_n)$.
For nonnegative integers $s_1,\ldots,s_n$, define
$s_+=\sum_{i=1}^ns_i$ and $a_+=\sum_{i=1}^na_i$. Then
\[
 \E\!\left[\prod_{i=1}^nX_i^{s_i}\right]
 =\frac{\prod_{i=1}^n(a_i)_{s_i}}{(a_+)_{s_+}},
\]
where $(q)_m=q(q+1)\cdots(q+m-1)$ denotes the Pochhammer symbol and
$(q)_0=1$.
\end{lemma}

\begin{lemma}[Application of the Dirichlet Mixed-Moment Formula]\label{lem:moment_application}
Let $\bm X\sim\Dir_n(a_1,\ldots,a_n)$ and $a_+=\sum_i a_i$.
Let $\bm r=(r_1,\ldots,r_n)$ and $\bm m=(m_1,\ldots,m_{n-1})$ have
nonnegative integer entries, and define $r_+=\sum_i r_i$ and
$m_+=\sum_{k=1}^{n-1}m_k$. Then
\[
 \E\!\left[\prod_{k=1}^{n-1}X_k^{r_k+m_k}X_n^{r_n}\right]
 =\frac{\prod_{k=1}^{n-1}(a_k)_{r_k+m_k}(a_n)_{r_n}}
 {(a_+)_{r_++m_+}}.
\]
\end{lemma}
\begin{proof}
By Lemma~\ref{lem:dirichlet_moments}, the Dirichlet mixed-moment formula
applies with powers $s_k=r_k+m_k$ for $k=1,\ldots,n-1$ and
$s_n=r_n$. Therefore,
\[
 s_+=\sum_{k=1}^{n-1}(r_k+m_k)+r_n
 =\sum_{i=1}^nr_i+\sum_{k=1}^{n-1}m_k=r_++m_+.
\]
Substituting these exponents into the Dirichlet mixed-moment formula gives
the result.
\end{proof}

\begin{lemma}[Multinomial Expansion]\label{lem:multinomial}
Let $p\in\mathbb N$ and let $y_1,\ldots,y_p$ be real or complex numbers.
For any nonnegative integer $m$,
\[
 \left(\sum_{k=1}^py_k\right)^m
 =\sum_{\substack{m_1,\ldots,m_p\geq0\\m_1+\cdots+m_p=m}}
 \frac{m!}{m_1!\cdots m_p!}\prod_{k=1}^py_k^{m_k}.
\]
\end{lemma}

\begin{lemma}[Carlson-Lauricella Relationship; \cite{Carlson1963, Carlson1998}]\label{lem:carlson}
Let $b_i>0$, $z_i>0$, $c=\sum_{i=1}^nb_i$ and $0\leq a<c$. Then
\[
 z_n^{-a}\FD^{(n-1)}\!\left(a;b_1,\ldots,b_{n-1};c;
 1-\frac{z_1}{z_n},\ldots,1-\frac{z_{n-1}}{z_n}\right)
 =R_{-a}(b_1,\ldots,b_n;z_1,\ldots,z_n),
\]
where the Lauricella function is evaluated by its defining series or its
real Euler continuation, and Carlson's multivariate $R$-function is
\[
 R_{-a}(\bm b;\bm z)
 =\frac1{B(\bm b)}\int_{S^{n-1}}
 \left(\sum_{i=1}^n z_iu_i\right)^{-a}
 \prod_{i=1}^n u_i^{b_i-1}\dd\bm u.
\]
Here $B(\bm b)=\prod_i\Gamma(b_i)/\Gamma(\sum_i b_i)$,
$S^{n-1}=\{\bm u:u_i>0,\ \sum_i u_i=1\}$, and
$\dd\bm u=\dd u_1\cdots\dd u_{n-1}$. For $a=0$, both sides equal one.
\end{lemma}

\begin{definition}[Lauricella Function of Type $D$; \cite{Carlson1963} and \cite{MatsumotoOhara2009}]\label{def:Lauricella}
Let $a,b_1,\ldots,b_p,c\in\mathbb C$, with
$c\notin\{0,-1,-2,\ldots\}$. The Lauricella function of type $D$ is
defined by
\[
 \FD^{(p)}(a;b_1,\ldots,b_p;c;x_1,\ldots,x_p)
 =\sum_{m_1,\ldots,m_p\geq0}
 \frac{(a)_{m_+}(b_1)_{m_1}\cdots(b_p)_{m_p}}
 {(c)_{m_+}m_1!\cdots m_p!}x_1^{m_1}\cdots x_p^{m_p},
\]
where $m_+=m_1+\cdots+m_p$. The series converges absolutely for
$\max_k|x_k|<1$. For real $c>a>0$ and real $x_k<1$, the Euler integral
\[
 \FD^{(p)}(a;b_1,\ldots,b_p;c;x_1,\ldots,x_p)
 =\frac{\Gamma(c)}{\Gamma(a)\Gamma(c-a)}
 \int_0^1t^{a-1}(1-t)^{c-a-1}\prod_{k=1}^p(1-x_kt)^{-b_k}\dd t
\]
gives the real continuation used in this paper. For $a=0$, the value is
one. In the cases $p=1$ and $p=2$, the function reduces to
${}_2F_1$ and Appell's $F_1$, respectively; the latter is also defined in
NIST DLMF, Section~16.13 in \cite{DLMF}.
\end{definition}

\bibliographystyle{unsrtnat}

\bibliography{main}
\end{document}